\documentclass{article}
\usepackage{mathrsfs}
\usepackage{amsfonts}
\usepackage{latexsym,amsfonts,amssymb,amsmath,amsthm,amscd,amsmath,graphicx}
\newtheorem{thm}{Theorem}[section]

\newtheorem{lem}[thm]{Lemma}

\theoremstyle{definition}

\theoremstyle{remark}
\newtheorem{rem}[thm]{Remark}
\numberwithin{equation}{section}

\newcommand{\p}{\partial}

\newcommand{\sV}{{\mathcal V}}
\newcommand{\sW}{{\mathcal W}}
\newcommand{\sF}{{\mathcal F}}
\newcommand{\Z}{{\mathbb Z}}
\newcommand{\N}{{\mathbb N}}
\newcommand{\x}{\xrightarrow}
\numberwithin{equation}{section}

\begin {document}
\topmargin= -.2in \baselineskip=20pt

\title{Calculation of local Fourier transforms for formal connections}
\author{Jiangxue Fang\\Chern Institute of Mathematics, Nankai University, Tianjin 300071,
P. R. China\\fangjiangx@yahoo.com.cn}
\date{}

\maketitle
\noindent\textbf{Abstract}\; We calculate the local Fourier
transforms for formal connections. In particular, we verify an
analogous conjecture of Laumon and Malgrange (\cite{La} 2.6.3).

\section{Introduction}
Let $k$ be an algebraic closed field of characteristic zero and let
$k((t))$ be the field of formal Laurent series in the variable $t$.
A formal connection on $k((t))$ is a pair $(M,t\p_t)$ consisting of
a finite dimensional $k((t))$-vector space $M$ and a $k$-linear map
$t\p_t:M\to M$ satisfying
$$t\p_t(fm)=t\p_t(f)m+f t\p_t(m)$$
for any $f\in k((t))$ and $m\in M$. In \cite{BE}, S. Bloch and H.
Esnault define local Fourier transforms $\sF^{(0,\infty)},\;
\sF^{(\infty, 0)},\;\sF^{(\infty,\infty)}$ for formal connections,
by analogy with the $\ell$-adic local Fourier transform considered
in \cite{La}. In \cite{La}, 2.6.3, Laumon and Malgrange give
conjectural formulas of local Fourier transforms for a class of
${\rm Q_\ell}$-sheaf. This results are proved by Lei Fu (\cite{F}).
In this paper, we prove an analogous conjecture of local Fourier
transform for formal connections. Actually, we can calculate local
Fourier transforms for any formal connections.

A key technical tool for the definitions of local Fourier transforms
of formal connections is the notion of good lattices pairs. By
definition in \cite{De}, Lemma 6.21, \emph{a pair of good lattices}
$\sV,\;\sW$ of $M$ is a pair of lattices in $M$ satisfying the
following conditions

(1) $\sV\subset\sW\subset M$

(2) $t\p_t(\sV)\subset\sW$

(3) For any $k\in\N$, the natural inclusion of complexes
$$(\sV\x{t\p_t}\sW)\to(\frac{1}{t^k}\sV\x{t\p_t}\frac{1}{t^k}\sW)$$
is a quasi-isomorphism.

Good lattices pairs $\sV,\;\sW$ exist. The number dim$_k\sW/\sV$ is
independent of the choice of good lattices pairs of $M$, and is
called the irregularity of $M$.

For any $f\in k((t)),$ denote by $[f]$ the formal connection on
$k((t))$ consisting of a one dimensional   $k((t))$-vector space
with a basis $e$ and a $k$-linear map $t\p_t:k((t)) e\to k((t)) e$
satisfying
$$t\p_t(ge)=(t\p_t(g)+fg)e$$ for any $g\in k((t)).$ Two such objects
$[f]$ and $[f']$ are isomorphic if and only if $f-f'\in tk[[t]]+\Z$.
Therefore the non-negative integer
$${\rm max}(0,-{\rm ord}_t(f))$$
is a well-defined invariant of the isomorphic class of $[f]$, and is
called the slope of $[f]$. Let $p$ be the slope of $[f]$. One can
verify $k[[t]] e,\;t^{-p}k[[t]] e$ is a good lattices pair of $[f]$.
So the irregularity coincides with the slope for any one dimensional
formal connection. The definition of slopes for arbitrary formal
connections is given in \cite{Ka-1}, (2.2.5). The irregularity of a
formal connection coincide with the sum of its slopes. Any formal
connection has a unique slope decomposition. So the slope of an
irreducible formal connection is equal to its irregularity divided
by its dimension. A formal connection is called regular if the
irregularity of this connection is equal to 0.

Throughout this paper, $r$ and $s$ are to be positive integers. Let
$t'$ be the Fourier transform coordinate of $t$. Write
$z=\frac{1}{t}$ and $z'=\frac{1}{t'}.$ Let
$$[r]:k((t))\hookrightarrow k((\sqrt[r]t))$$ be the natural
inclusion of fields. Let $T=\sqrt[r]t$ and let $\alpha$ be a formal
Laurent series in $k((T))$ of order $-s$ with respect to $T$. Let
$R$ be a regular formal connection on $k((T))$. In this paper, we
calculate the local Fourier transform
$$\sF^{(0,\infty)}\Big([r]_*\Big([T\p_T(\alpha)]\otimes_{k((T))} R\Big)\Big).$$
Similarly, let $k((z))$ be the field of formal Laurent series in the
variable $z.$ Let $$[r]:k((z))\hookrightarrow
k((\frac{1}{\sqrt[r]t}))$$ be the natural inclusion of fields. Let
$Z=\frac{1}{\sqrt[r]t}$ and let $\alpha$ be a formal Laurent series
in $k((Z))$ of order $-s$ with respect to $Z$. Let $R$ be a regular
formal connection on $k((Z))$. We also calculate the local Fourier
transforms
\begin{eqnarray*}
&&\sF^{(\infty,0)}\Big([r]_*\Big([Z\p_Z(\alpha)]\otimes_{k((Z))}R\Big)\Big)\hbox{ if }r>s;\\
&&\sF^{(\infty,\infty)}\Big([r]_*\Big([Z\p_Z(\alpha)]\otimes_{k((Z))}R\Big)\Big)\hbox{
if }r<s.
\end{eqnarray*}

 We refer the reader to \cite{BE} for
the definitions and properties of local Fourier transforms. The main
results of this paper are the following three theorems.
\bigskip

\noindent\textbf{Theorem 1.} \emph{Given a formal Laurent series
$\alpha$ in $k((\sqrt[r]t))$ of order $-s$ with respect to
$\sqrt[r]t$, consider the following system of equations
\begin{eqnarray}\label{1}\left\{\begin{array}{ll}
\p_t(\alpha(\sqrt[r]t))+t'=0,\\
\alpha(\sqrt[r]t)+tt'=\beta(\frac{1}{\sqrt[r+s]{t'}}).
\end{array}\right.\end{eqnarray}
Using the first equation, we find an expression of $\sqrt[r]t$ in
terms of $\frac{1}{\sqrt[r+s]{t'}}$. We then substitute this
expression into the second equation to get
$\beta(\frac{1}{\sqrt[r+s]{t'}}),$ which is a formal Laurent series
in $k((\frac{1}{\sqrt[r+s]{t'}}))$ of order $-s$ with respect to
$\frac{1}{\sqrt[r+s]{t'}}$. Let $T=\sqrt[r]t$ and let
$Z'=\frac{1}{\sqrt[r+s]{t'}}$. For any regular formal connection $R$
on $k((T))$, we have
$$\sF^{(0,\infty)}\Big([r]_*\Big([T\p_T(\alpha)]\otimes_{k((T))} R\Big)\Big)=[r+s]_*\Big([Z'\p_{Z'}(\beta)+\frac{s}{2}]\otimes_{K((Z'))} R\Big),$$
where the right $R$ means the formal connection on $k((Z'))$ after
replacing the variable $T$ with $Z'$.}

\noindent\textbf{Theorem 2.} \emph{Suppose $r>s$. Given a formal
Laurent series $\alpha$ in $k((\frac{1}{\sqrt[r]t}))$ of order $-s$
with respect to $\frac{1}{\sqrt[r]t}$, consider the following system
of equations
\begin{eqnarray}\label{2}\left\{\begin{array}{ll}
\p_t(\alpha(\frac{1}{\sqrt[r]t}))+t'=0,\\
\alpha(\frac{1}{\sqrt[r]t})+tt'=\beta(\sqrt[r-s]{t'}).
\end{array}\right.\end{eqnarray}
Using the first equation, we find an expression of
$\frac{1}{\sqrt[r]t}$ in terms of $\sqrt[r-s]{t'}.$ We then
substitute this expression into the second equation to get
$\beta(\sqrt[r-s]{t'}),$ which is formal Laurent series in
$k((\sqrt[r-s]{t'}))$ of order $-s$ with respect to
$\sqrt[r-s]{t'}$. Let $Z=\frac{1}{\sqrt[r]t}$ and let
$T'=\sqrt[r-s]{t'}.$ For any regular formal connection $R$ on
$k((Z))$, we have
$$\sF^{(\infty,0)}\Big([r]_*\Big([Z\p_Z(\alpha)]\otimes_{k((Z))} R \Big)\Big)=[r-s]_*\Big([T'\p_{T'}(\beta)+\frac{s}{2}]\otimes_{k((T'))} R\Big),$$
where the right $R$ means the formal connection on $k((T'))$ after
replacing the variable $Z$ with $T'$. }

\noindent\textbf{Theorem 3.} \emph{Suppose $r<s$. Given a formal
Laurent series $\alpha$ in $k((\frac{1}{\sqrt[r]t}))$ of order $-s$
with respect to $\frac{1}{\sqrt[r]t}$, consider the following system
of equations
\begin{eqnarray}\label{2.}\left\{\begin{array}{ll}
\p_t(\alpha(\frac{1}{\sqrt[r]t}))+t'=0,\\
\alpha(\frac{1}{\sqrt[r]t})+tt'=\beta(\frac{1}{\sqrt[s-r]{t'}}).
\end{array}\right.\end{eqnarray}
Using the first equation, we find an expression of
$\frac{1}{\sqrt[r]t}$ in terms of $\frac{1}{\sqrt[s-r]{t'}}.$ We
then substitute this expression into the second equation to get
$\beta(\frac{1}{\sqrt[s-r]{t'}}),$ which is a formal Laurent series
in $k((\frac{1}{\sqrt[s-r]{t'}}))$ of order $-s$ with respect to
$\frac{1}{\sqrt[s-r]{t'}}$. Let $Z=\frac{1}{\sqrt[r]t}$ and let
$Z'=\frac{1}{\sqrt[s-r]{t'}}.$ For any regular formal connection $R$
on $k((Z))$, we have
$$\sF^{(\infty,\infty)}\Big([r]_*\Big([Z\p_Z(\alpha)]\otimes_{k((Z))} R\Big)\Big)=[s-r]_*\Big([Z'\p_{Z'}(\beta)+\frac{s}{2}]\otimes_{k((Z'))} R\Big),$$
where the right $R$ means the formal connection on $k((Z'))$ after
replacing the variable $Z$ with $Z'$.}

When $R$ is trivial, the above three theorems are conjectured by
Laumon and Malgrange (\cite{La} 2.6.3) except the term $\frac{s}{2}$
is missing in the conjecture. Any formal connection on $k((t))$ is a
direct sum of indecomposable connections. As in \cite{BBE}, section
5.9, any indecomposable connection $M=N\otimes R$, where $R$ is
regular and $N=[d]_*L$ where $L$ is a one dimensional connection on
a finite extension $[d]: k((t))\to k((t^{\frac{1}{d}}))$. So we can
calculate local Fourier transform for all formal connections.

\textbf{\;\;Acknowledgements.} It is a great pleasure to thank my
advisor Lei Fu for his guidance and support during my graduate
studies. In \cite{Sa}, Claude Sabbah proves these results of local
Fourier transforms for formal connections with a geometric method.
Our method is elementary and directly.

\section{Proofs of Theorems 1, 2}
Given a formal Laurent series $\alpha$ in the variable $\sqrt[r]t$
of order $-s$, consider the system of equations (\ref{1}). We
express $\sqrt[r]t$ as a formal Laurent series in
$\frac{1}{\sqrt[r+s]{t'}}$ of order 1 using the first equation and
then substitute this expression into the second equation to get
$\beta\in k((\frac{1}{\sqrt[r+s]{t'}}))$. We have
\begin{eqnarray}\label{*}
\p_{t'}(\beta)&=&\p_{t'}\big(\alpha(\sqrt[r]t)+tt'\big)
=\p_t\big(\alpha(\sqrt[r]t)\big)\frac{dt}{dt'}+t'\frac{dt}{dt'}+t\\\nonumber&=&\Big(\p_t\big(\alpha(\sqrt[r]t)\big)+t'\Big)\frac{dt}{dt'}+t=t.
\end{eqnarray}
It follows that  $\beta$ is a formal Laurent series in
$\frac{1}{\sqrt[r+s]{t'}}$ of order $-s$. Let $T=\sqrt[r]t$ and
$Z'=\frac{1}{\sqrt[r+s]{t'}}$. Set
$$a(T)=-T^st\p_t(\alpha) \hbox{ and } b(Z')=Z'^st'\p_{t'}(\beta).$$
Then $a(T)$ is a formal power series in $T$ of order 0 and $b(Z')$
is a formal power series in $Z'$ of order 0. From the system of
equations (\ref{1}) and (\ref{*}), we get
\begin{eqnarray}\label{21}\left\{\begin{array}{ll}
a(T)=(\frac{T}{Z'})^{r+s}\\
b(Z')=(\frac{T}{Z'})^r.
\end{array}\right.\end{eqnarray}
To prove Theorem 1, it suffices to prove the following theorem.

\smallskip

\noindent\textbf{Theorem 1$'$.} \emph{Given a formal power series
$a(T)=\sum_{i\geq0}a_iT^i$ with $a_i\in k$ and $a_0\neq0$, solve the
system of equations (\ref{21}) to get $b(Z')=\sum_{i\geq0}b_iZ'^i$
for some $b_i\in k$. Then $b_s=\frac{r}{r+s}a_s$ and
\begin{eqnarray*}
&&\sF^{(0,\infty)}\Big([r]_*[-r(a_0T^{-s}+a_1T^{1-s}+\ldots+a_s)]\Big)\\
&=&[r+s]_*[-(r+s)(b_0Z'^{-s}+b_1Z'^{1-s}+\ldots+b_s)+\frac{s}{2}].
\end{eqnarray*}}
In fact, suppose Theorem 1$'$ holds. Let $c$ be an element in $k.$
By remark \ref{r} we shall prove later, for
$a(T)=-T^st\p_t(\alpha)-\frac{c}{r}T^s$, we can get a solution
$b(Z')$ of the system of equations (\ref{21}) such that
$$b(Z')\equiv Z'^st'\p_{t'}(\beta)-\frac{c}{r+s}Z'^s \hbox{ mod. }Z'^{s+1}.$$
 Then
\begin{eqnarray*}
&&\sF^{(0,\infty)}\Big([r]_*[T\p_T(\alpha)+c\Big)\\
&=&\sF^{(0,\infty)}\Big([r]_*[-rT^{-s}(-T^st\p_t(\alpha)-\frac{c}{r}T^s)]\Big)\\
&=&[r+s]_*[-(r+s)Z'^{-s}(Z'^st'\p_{t'}(\beta)-\frac{c}{r+s}Z'^s)]\\
&=&[r+s]_*[Z'\p_{Z'}(\beta)+c].
\end{eqnarray*}
So Theorem 1 holds for $R=[c]$. As in \cite{BBE}, section 5.9, every
irreducible regular formal connection $N$ on $k((T))$ is $[d]_*L$,
where $L$ is a one dimensional formal connection on a finite
extension $[d]: k((T))\to k((T^{\frac{1}{d}}).$ So $L$ is regular,
we have $L=[c]$ for some $c\in k$. Then $N=[d]_*[c]=\oplus_{1\leq
i\leq d}[c+\frac{i}{d}].$ We have $d=1$ because $N$ is irreducible.
This shows that every irreducible regular formal connection is
isomorphic to the one dimensional connection $[c]$ for some $c\in
k.$ So every regular formal connection is a successive extension of
connections of the type $[c]$. Since $\sF^{(0,\infty)}$ is
functoriel and exact, Theorem 1 holds for any regular formal
connection $R$ on $k((T))$.

\begin{rem}\label{rem}
If $a_s=0$, then there exists $\alpha\in k((\sqrt[r]t))$ such that
$a(T)=-T^st\p_t(\alpha).$ Using the first equation of (\ref{21}), we
find an expression of $T$ in terms of $Z'$. We then substitute this
expression into the second equation of (\ref{21}) to get $b(Z').$
This expression also satisfies the first equation of (\ref{1}). We
then substitute this expression into the second equation of
(\ref{1}) to get $\beta(Z').$ By (\ref{*}), we have
$$b(Z')=\sum_{i\geq0}b_iZ'^i=Z'^st'\p_{t'}(\beta).$$
This shows $b_s=0$.
\end{rem}

\begin{rem}\label{r}
Solving the first equation of (\ref{21}), we get
$T=\sum_{i\geq0}\lambda_iZ'^{i+1}$ with $\lambda_0=\sqrt[r+s]{a_0}.$
The solution is not unique and different solutions differ by an
$r+s$-th root of unity. As long as $\lambda_0$ is chosen to be an
$r+s$-th root of $a_0$, for each $i$, $\lambda_i$ depends only on
$a_0,\ldots,a_i$. We have $b(Z')=(\sum_{i\geq0}\lambda_iZ'^i)^r,$
and for each $i$, $b_i$ depends only on
$\lambda_0,\ldots,\lambda_i$. Therefore as long as we fix an
$r+s$-th root of $a_0$, for each $i$, $b_i$ depends only on
$a_0,\ldots,a_i$. So to prove Theorem $1'$, we can assume
$a(T)=\sum_{0\leq i\leq s}a_iT^i.$
\end{rem}

\begin{rem}\label{re}
Solving the first equation of (\ref{21}), we get
$T=\sum_{i\geq0}\lambda_iZ'^{i+1}$ for some $\lambda_j\in k.$ Then
$\lambda_0$ is an $r+s$-th root of $a_0.$ Then
$\sum_{i\geq0}b_iZ'^i=(\sum_{i\geq0}\lambda_iZ'^i)^r$. Choose
$a_0',\ldots,a_s'\in k$ such that $a'_i=a_i$ for all $0\leq i<s$ and
$a_s'=0$. For $a(T_1)=\sum_{0\leq i\leq s}a_i'T_1^i,$ consider the
system of equations (\ref{21}) if the variable $T$ is changed by
$T_1$. Using the first equation, we can express $T_1$ as
$\sum_{i\geq0}\lambda_i'Z'^{i+1}$ with $\lambda'_0=\lambda_0$. Then
we have $\sum_{i\geq0}b_i'Z'^i=(\sum_{i\geq0}\lambda_i'Z'^i)^r.$
Remark \ref{rem} shows $b'_s=0$. Since $a_i=a_i'$ for $0\leq i<s$,
we have $\lambda_i=\lambda_i'$ for all $0\leq i<s$. That is,
$$T\equiv T_1\hbox{ mod. }Z'^{s+1} \hbox{ and } T\equiv T_1\equiv \lambda_0Z'\hbox{ mod. } Z'^2.$$
Comparing coefficients of $Z'^s$ on both sides of
$$\sum_{i\geq0}a_iT^i=\Big(\sum_{i\geq0}\lambda_iZ'^i\Big)^{r+s}\hbox{ and }
\sum_{0\leq i\leq
s}a'_iT_1^i=\Big(\sum_{i\geq0}\lambda'_iZ'^i\Big)^{r+s},$$ we have
$$a_s\lambda_0^s=(a_s-a_s')\lambda_0^s=(r+s)(\lambda_s-\lambda_s')\lambda_0^{r+s-1}.$$
Comparing coefficients of $Z'^s$ on both sides of
$$\sum_{i\geq0}b_iZ'^i=\Big(\sum_{i\geq0}\lambda_iZ'^i\Big)^r\hbox{ and
}\sum_{i\geq0}b'_iZ'^i=\Big(\sum_{i\geq0}\lambda'_iZ'^i\Big)^r,$$ we
have
$$b_s=b_s-b_s'=r(\lambda_s-\lambda_s')\lambda_0^{r-1}.$$
This proves $b_s=\frac{r}{r+s}a_s.$
\end{rem}

\begin{rem}\label{r2.1} Set $f=a_0T^{-s}+a_1T^{1-s}+\ldots+a_s$. Let
$$H=\{\sigma\in {\rm Gal}\big(k((T))/k((t))\big)|\sigma(f)=f\}.$$
We call $f$ is irreducible with respect to the Galois extension
$k((T))/k((t))$ if $\#H=1$. Then $f$ is irreducible if and only if
the connection $[r]_*[-rf]$ is irreducible.
\end{rem}

\begin{lem}\label{l52}
If Theorem $1'$ holds for irreducible $f$, then it holds for all f.
\end{lem}

\begin{proof}
By Remark \ref{r}, we can assume $a(T)=\sum_{0\leq i\leq s}a_iT^i.$
Keep the notation in Remark \ref{r2.1}. Set $p=\#H$. Then $p|r$. Let
$\eta$ be a primitive $r$-th root of unity in $k$. Then
$a_i\eta^{\frac{r}{p}(i-s)}=a_i$ for all $0\leq i\leq s$. So $a_i=0$
or $p|i-s$. In particular, $p|s$ since $a_0\neq0$. Let $\tau=T^p$
and $\tau'=Z'^p$. Then
$$f=a_0\tau^{-\frac{s}{p}}+a_p\tau^{1-\frac{s}{p}}+\ldots+a_s$$ and it
is irreducible with respect to the Galois extension
$k((\tau))/k((t)).$ For $a(\tau)=\sum_{0\leq i\leq
\frac{s}{p}}a_{pi}\tau^i$, suppose
$b(\tau')=\sum_{i\geq0}b_{pi}\tau'^i$ is a solution of the following
system of equation
\begin{eqnarray}\label{tau}\left\{\begin{array}{ll}
a(\tau)=(\frac{\tau}{\tau'})^{\frac{r+s}{p}}\\
b(\tau')=(\frac{\tau}{\tau'})^{\frac{r}{p}}.
\end{array}\right.\end{eqnarray}
Then $b_s=\frac{r}{r+s}a_s$ and $b(Z')=\sum_{i\geq0}b_{pi}Z'^{pi}$
is a solution of the system of equations (\ref{21}). For
$a(\tau)=\sum_{0\leq i\leq
\frac{s}{p}}a_{pi}\tau^i-\frac{j}{r}\tau^{\frac{s}{p}}\;(1\leq j\leq
p)$, by Remark \ref{r} and \ref{re}, we can find a solution
$b(\tau')$ of the system of equations (\ref{tau}) such that
$$b(\tau')\equiv \sum_{0\leq i\leq
\frac{s}{p}}b_{pi}\tau'^i-\frac{j}{r+s}\tau'^{\frac{s}{p}}\hbox{
mod. }\tau'^{\frac{s}{p}+1}.$$ Applying Theorem 1$'$ to the system
of equations (\ref{tau}) for $a(\tau)=\sum_{0\leq i\leq
\frac{s}{p}}a_{pi}\tau^i-\frac{j}{r}\tau^{\frac{s}{p}}\;(1\leq j\leq
p)$, we have
\begin{eqnarray*}
&&\sF^{(0,\infty)}\Big([r]_*[-rf]\Big)\\
&=&\sF^{(0,\infty)}\Big([\frac{r}{p}]_*[p]_*[-r(a_0T^{-s}+a_pT^{p-s}+\ldots+a_s)]\Big)\\
&=&\bigoplus_{1\leq j\leq
p}\sF^{(0,\infty)}\Big([\frac{r}{p}]_*[-\frac{r}{p}(a_0\tau^{-\frac{s}{p}}+a_p\tau^{\frac{p-s}{p}}+\ldots+a_s)+\frac{j}{p}]\Big)\\
&=&\bigoplus_{1\leq j\leq
p}[\frac{r+s}{p}]_*[-\frac{r+s}{p}(b_0\tau'^{-\frac{s}{p}}+b_p\tau'^{\frac{p-s}{p}}+\ldots+b_s)+\frac{j}{p}+\frac{s}{2p}]\\
&=&[\frac{r+s}{p}]_*[p]_*[-(r+s)(b_0Z'^{-s}+b_pZ'^{p-s}+\ldots+b_s)+\frac{s}{2}]\\
&=&[r+s]_*[-(r+s)(b_0Z'^{-s}+b_pZ'^{p-s}+\ldots+b_s)+\frac{s}{2}].
\end{eqnarray*}
\end{proof}
From now on, we assume $f$ is irreducible.

\medskip

Let's describe the connection $\sF^{(0,\infty)}\big([r]_*[-rf]\big)$
on $k((z'))$.

The formal connection $[-rf]$ on $k((T))$ consist of a one
dimensional $k((T))$-vector space with a basis $e$ and a $k$-linear
map $T\p_T:k((T))e\to k((T))e$ satisfying
$$T\p_T(ge)=(T\p_T(g)-rfg)e$$
for any $g\in k((T)).$ Since the formal connection $[-rf]$ on
$k((T))$ has slope $s$, we get $k[[T]] e,\;T^{-s}k[[T]] e$ is a good
lattices pair for it. Identify $[r]_*[-rf]$ with $k((T)) e$ as
$k((t))$-vector spaces. Then the formal connection $[r]_*[-rf]$ has
pure slope $\frac{s}{r}$ and $k[[T]]e,\;T^{-s}k[[T]] e$ is a good
lattices pair for this connection. The action of the differential
operator $t\p_t$ on $k((T))e$ is given by
$$t\p_t(g e)=(t\p_t(g)-fg)e$$
for any $g\in k((T)).$ So we have
\begin{eqnarray*}
&&({\p_t}\circ
t)(T^{-i}e)=\frac{r-i}{r}T^{-i}e-(a_0T^{-(s+i)}e+\ldots+a_sT^{-i}e)\;(1\leq i\leq r),\\
&&t\cdot T^{-i}e=T^{-(i-r)}e\;(r+1\leq i\leq r+s).
\end{eqnarray*}
By \cite{BE}, Proposition 3.7, the map
$$\iota:k((T))e\to\sF^{(0,\infty)}\Big([r]_*[-rf]\Big)$$ is an isomorphism
of $k$-vector spaces. By \cite{BE}, Lemma 2.4, $(\iota
T^{-1}e,\ldots,\iota T^{-(r+s)}e)$ is a basis of
$\sF^{(0,\infty)}\Big([r]_*[-rf]\Big)$ over $k((z'))$. Then by the
relation $\iota\circ t=-z'^2\p_{z'}\circ\iota$ and
$\iota\circ\p_t=-\frac{1}{z'}\circ\iota$ in \cite{BE}, Proposition
3.7, the matrix of the connection
$\sF^{(0,\infty)}\big([r]_*[-rf]\big)$ with respect to the
differential operator $z'\p_{z'}$ and the basis $(\iota
T^{-1}e,\ldots,\iota T^{-(r+s)}e)$ is \[ -\bordermatrix{
&\overbrace{\ \ \ \ \ \ \ \ \ \ \ \ \ }^{r} &\overbrace{\ \ \ \ \ \
\ \ \ \ \ \ \ }^{s}&\cr
&\begin{array}{ccc}a_s&&\\a_{s-1}&\ddots&\\\vdots&\ddots&a_s\end{array}&\begin{array}{ccc}
\frac{1}{z'}\\&\ddots\\&&\frac{1}{z'}\end{array}\cr
&\begin{array}{ccc}
a_0&&a_{s-1}\\&\ddots&\vdots\\&&a_0\end{array}&\begin{array}{ccc}
 \\&\ \\&& \end{array}\cr } +{\rm
diag}\{\frac{r-1}{r},\ldots,\frac{1}{r},0,\ldots,0\}.\] Then the
matrix of the connection
$$[r+s]^*\Big(\sF^{(0,\infty)}\Big([r]_*[-rf]\Big)\Big)=k((Z'))\otimes_{k((z'))}\sF^{(0,\infty)}\Big([r]_*[-rf]\Big)$$
with respect to the differential operator $Z'\p_{Z'}$ and the basis
$(Z'\otimes\iota T^{-1}e,\ldots,Z'^{r+s}\otimes\iota T^{-(r+s)}e)$
is\begin{eqnarray*} &&-\frac{r+s}{Z'^s}\left(\begin{array}{cccccc}
a_sZ'^s&&& 1\\
a_{s-1}Z'^{s-1}&\ddots&&&\ddots\\
\vdots&\ddots&a_sZ'^s&&&1\\
a_0&\ddots&a_{s-1}Z'^{s-1}&&\\
&\ddots&\vdots&&\\
&&a_0
\end{array}\right)\\
&&+(r+s){\rm diag}\{\frac{r-1}{r},\ldots,\frac{1}{r},0,\ldots,0\}
+{\rm diag}\{1, \ldots, r+s \}.
\end{eqnarray*}
We can write this matrix as $(r+s)B-(r+s)\sum_{0\leq i\leq
s}Z'^{i-s}A_{i}$ for some matrices $A_i$ and $B$ with entries in
$k$, where
\begin{eqnarray*}
A_0&=&\left(\begin{array}{cc} 0&I_s\\a_0I_r&0\end{array}\right),\\
 B&=&{\rm diag}\{\frac{r-1}{r},\ldots,\frac{1}{r},0,\ldots,0\}+\frac{1}{r+s}{\rm
diag}\{1,\ldots,r+s\}.
\end{eqnarray*}
Let $V$ be the $k$-vector subspace of
$[r+s]^*\Big(\sF^{(0,\infty)}\Big([r]_*[-rf]\Big)\Big)$ generated by
$Z'^i\otimes\iota T^{-i}e$ $(1\leq i\leq r+s)$. With respect to this
basis, $V$ can be identified with the $k$-vector space of column
vectors in $k$ of length $r+s$. The action of the differential
operator $Z'\p_{Z'}$ on elements of $V$ can be written as
$$Z'\p_{Z'}(v)=(r+s)B(v)-(r+s)\sum_{0\leq i\leq s}Z'^{i-s}A_i(v).$$

\begin{lem}\label{l53}
Suppose $f$ is irreducible in the sense of Remark \ref{r2.1}. Given
$\alpha_0,\ldots,\alpha_s\in k$, the following three conditions are
equivalent:

$(1)$
$\sF^{(0,\infty)}\big([r]_*[-rf]\big)=[r+s]_*[-(r+s)\sum_{0\leq
i\leq s}\alpha_iZ'^{i-s}].$

$(2)$ $[-(r+s)\sum_{0\leq i\leq s}\alpha_iZ'^{i-s}]$ is a
subconnection of
$[r+s]^*\Big(\sF^{(0,\infty)}\Big([r]_*[-rf]\Big)\Big)$.

$(3)$ There exist an integer $N$ and $v_0,\ldots,v_s\in V$ such that
$v_0\neq0$ and
\begin{eqnarray}\label{5.3}\left\{\begin{array}{ll}
\sum_{0\leq i\leq k}(A_i-\alpha_i)v_{k-i}=0\;(0\leq k\leq s-1);\\
\sum_{0\leq i\leq
s-1}(A_i-\alpha_i)v_{s-i}+(A_s-B-\alpha_s-\frac{N}{r+s})v_0=0.
\end{array}\right.\end{eqnarray}
\end{lem}

\begin{proof}
Since $f$ is irreducible, the connection $[r]_*[-rf]$ on $k((t))$ is
irreducible with pure slope $\frac{s}{r}$. By \cite{BE}, Proposition
3.14, the connection $\sF^{(0,\infty)}\Big([r]_*[-rf]\Big)$ on
$k((z'))$ is irreducible with pure slope $\frac{s}{r+s}$. As in the
proof of \cite{BE}, Lemma 3.3, we have
$$\sF^{(0,\infty)}\Big([r]_*[-rf]\Big)=[r+s]_*[-(r+s)\sum_{0\leq i\leq s}\varrho_{i}Z'^{i-s}]$$
for some $\varrho_0,\ldots,\varrho_s\in k$ with $\varrho_0\neq0.$
Let $\mu$ be a primitive $(r+s)$-th root of unity in $k$. Then
\begin{eqnarray*}
&&[r+s]^*\Big(\sF^{(0,\infty)}\Big([r]_*[-rf]\Big)\Big)\\
&=&\bigoplus_{1\leq j\leq
r+s}[-(r+s)(\mu^{-js}\varrho_0Z'^{-s}+\mu^{j(1-s)}\varrho_1Z'^{1-s}+\ldots+\varrho_s)].
\end{eqnarray*}
So there are $r+s$ one dimensional  subconnections of
$[r+s]^*\Big(\sF^{(0,\infty)}\Big([r]_*[-rf]\Big)\Big)$ which are
not isomorphic to each other.

$(1)\Rightarrow(2)$ is trivial. For $(2)\Rightarrow(1)$, assume that
$[-(r+s)\sum_{0\leq i\leq s}\alpha_iZ'^{i-s}]$ is a subconnection of
$[r+s]^*\Big(\sF^{(0,\infty)}\Big([r]_*[-rf]\Big)\Big)$. Then
$$[-(r+s)\sum_{0\leq i\leq
s}\alpha_iZ'^{i-s}]=[-(r+s)\sum_{0\leq i\leq
s}\mu^{j(i-s)}\varrho_{i}Z'^{i-s}]$$ for some $1\leq j\leq r+s$.
Then
\begin{eqnarray*}
\sF^{(0,\infty)}\Big([r]_*[-rf]\Big)&=&[r+s]_*[-(r+s)\sum_{0\leq
i\leq
s}\mu^{j(i-s)}\varrho_{i}Z'^{i-s}]\\
&=&[r+s]_*[-(r+s)\sum_{0\leq i\leq s}\alpha_iZ'^{i-s}].
\end{eqnarray*}

For $(2)\Rightarrow(3)$, assume that $[-(r+s)\sum_{0\leq i\leq
s}\alpha_iZ'^{i-s}]$ is a subconnection of
$[r+s]^*\Big(\sF^{(0,\infty)}\Big([r]_*[-rf]\Big)\Big).$ This means
that there is a nonzero map of connections
$$\phi:[-(r+s)\sum_{0\leq i\leq s}\alpha_iZ'^{i-s}]\to
[r+s]^*\Big(\sF^{(0,\infty)}\Big([r]_*[-rf]\Big)\Big).$$ The
connection $[-(r+s)\sum_{0\leq i\leq s}\alpha_iZ'^{i-s}]$ consist of
a one dimensional $k((Z'))$-vector space with a basis $\varepsilon$
and a $k$-linear map $Z'\p_{Z'}:k((Z'))\varepsilon\to
k((Z'))\varepsilon$ satisfying
$$Z'\p_{Z'}(g\varepsilon)=\Big(Z'\p_{Z'}(g)-(r+s)g\sum_{0\leq i\leq s}\alpha_iZ'^{i-s}\Big)\varepsilon$$
for any $g\in k((Z')).$ Suppose $\phi(\varepsilon)=\sum_{0\leq
i}Z'^{i+N}v_i$ for some integer $N$ and some $v_i\in V$ with
$v_0\neq0$. Then
\begin{eqnarray*}
&&-(r+s)\sum_{0\leq i\leq s}\alpha_iZ'^{i-s}\sum_{0\leq i}Z'^{i+N}v_i=\phi\big(Z'\p_{Z'}(\varepsilon)\big)\\
&=&Z'\p_{Z'}\big(\phi(\varepsilon)\big)=Z'\p_{Z'}\big(\sum_{0\leq i}Z'^{i+N}v_i\big)\\
&=&\sum_{0\leq i}Z'^{i+N}\big((r+s)B+i+N\big)v_i-(r+s)\sum_{0\leq
i}Z'^{i+N}\sum_{0\leq j\leq s} Z'^{j-s}A_j(v_i).
\end{eqnarray*}
Comparing coefficients of $Z'^i$, for $N-s\leq i\leq N$ on each
side, we get the system of equations (\ref{5.3}). This proves
$(2)\Rightarrow(3)$. So for
$\alpha_0=\mu^{-sj}\varrho_0,\alpha_1=\mu^{(1-s)j}\varrho_1,\ldots,\alpha_s=\varrho_s$,
the system of equations (\ref{5.3}) holds for some $N\in\Z$ and some
$v_0,\ldots,v_s\in V$ with $v_0\neq0$. These $(s+1)$-tuples
$(\mu^{-sj}\varrho_0,\mu^{(1-s)j}\varrho_1,\ldots,\varrho_s)$
$(1\leq j\leq r+s)$ are pairwise distinct, since $f$ is irreducible.
Lemma \ref{l54} shows that there are at most $r+s$ $(s+1)$-tuples
$(\alpha_0,\ldots,\alpha_s)$ such that the system of equations
(\ref{5.3}) holds for $N=0$ and some $v_0,\ldots,v_s\in V$ with
$v_0\neq 0$. This proves $(3)\Rightarrow(2)$.
\end{proof}

\noindent\textbf{Hensel's lemma.}  \emph{Let $E$ be a finite
dimensional $k$-vector space. Suppose $D$ is a $k[[t]]$-linear
endomorphism of $E\otimes_kk[[t]]$. Write the action of $D$ on
elements of $E$:
$$D(v)=\sum_{i\geq0}t^iD_i(v), \hbox{ for unique elements }D_i\in {\rm
End}_k(E).$$ Suppose the characteristic polynomial of $D_0$ has a
simple root $\alpha_0$ in $k$. Then}

\emph{{\rm (1)} The equation $$(D-\alpha)(u)=0$$ has a solution
$\alpha\in k[[t]]$ with constant term $\alpha_0$ and $0\neq u\in
E\otimes_kk[[t]].$ In this case, $\alpha$ is uniquely determined by
$\alpha_0.$}

\emph{{\rm (2)} Let $k$ be a positive integer. The following systems
of equations
$$\sum_{0\leq i\leq j}(D_i-\alpha_i)u_{j-i}=0\;(0\leq j\leq k)$$ has
a solution $\alpha_1,\ldots,\alpha_k\in k$; $u_0,\ldots,u_k\in E$
with $u_0\neq0$. In this case, $\alpha_1,\ldots,\alpha_k$ are
uniquely determined by $\alpha_0$.}

\begin{proof}
The proof is similar to that of \cite{S}, Proposition 7, p. 34.
\end{proof}

\begin{lem}\label{l54}
Given $\alpha_0,\ldots,\alpha_s\in k$, there exist
$v_0,\ldots,v_s\in V$ such that $v_0\neq0$ and
\begin{eqnarray}\label{513}\left\{\begin{array}{ll}
\sum_{0\leq i\leq k}(A_i-\alpha_i)v_{k-i}=0\;(0\leq k\leq s-1),\\
\sum_{0\leq i\leq s-1}(A_i-\alpha_i)v_{s-i}+(A_s-B-\alpha_s)v_0=0
\end{array}\right.\end{eqnarray}
if and only if there exist $v_0',\ldots,v_s'\in V$ such that
$v_0'\neq0$ and
\begin{eqnarray}\label{54}\left\{\begin{array}{ll}
\sum_{0\leq i\leq k}(A_i-\alpha_i)v'_{k-i}=0\;(0\leq k\leq s-1),\\
\sum_{0\leq i\leq
s-1}(A_i-\alpha_i)v'_{s-i}+(A_s-\frac{2r+s}{2r+2s}-\alpha_s)v_0'=0.
\end{array}\right.\end{eqnarray}
Moreover, there are at most $r+s$ $(s+1)$-tuples
$(\alpha_0,\ldots,\alpha_s)$ in $k$ such that the system of
equations (\ref{513}) (resp. (\ref{54})) holds for some
$v_0,\ldots,v_s\in V$ with $v_0\neq0$ (resp. $v'_0,\ldots,v'_s\in V$
with $v_0'\neq0$).
\end{lem}

\begin{proof}
Let $\mu$ be a primitive $(r+s)$-th root of unity in $k$. We fix an
$(r+s)$-th root $a_0^{\frac{1}{r+s}}$ of $a_0.$ For any $1\leq j\leq
r+s,$ set $e_j$ to be the column vector
$(\mu^ja_0^{\frac{1}{r+s}},\ldots,\mu^{j(r+s-1)}a_0^{\frac{r+s-1}{r+s}},a_0)$
and $\varepsilon_j$ the row vector
$(\mu^{-j}a_0^{-\frac{1}{r+s}},\ldots,\mu^{-j(r+s-1)}a_0^{-\frac{r+s-1}{r+s}},a_0^{-1})$.
Then
$$A_0\cdot e_j=\mu^{rj}a_0^{\frac{r}{r+s}}\cdot e_j, \;
\varepsilon_j\cdot
A_0=\mu^{rj}a_0^{\frac{r}{r+s}}\cdot\varepsilon_j, \;
\varepsilon_i\cdot e_j=(r+s)\delta_{ij}.$$  Set $d=(r,s)$. We get
ker$(A_0-\mu^{rj}a_0^{\frac{r}{r+s}})$ is generated by those $e_k$
with $r+s|(k-j)d$, and im$(A_0-\mu^{rj}a_0^{\frac{r}{r+s}})$ is
generated by the other $e_k$'s. Then
$${\rm im}(A_0-\mu^{rj}a_0^{\frac{r}{r+s}})=\{v\in V|\varepsilon_k\cdot
v=0 \hbox { for all $k$ satisfying } r+s|(k-j)d\}.$$ For the only if
part, suppose the system of equations (\ref{513}) holds for some
$v_0,\ldots,v_s\in V$ with $v_0\neq0$. In particular,
$(A_0-\alpha_0)v_0=0$. Then $\alpha_0=\mu^{rj}a_0^{\frac{r}{r+s}}$
for some integer $j$ and then $v_0=\sum_{r+s|(i-j)d}\gamma_ie_i$ for
some $\gamma_i\in k$. For any $1\leq k,\;l\leq r+s$, we have
\begin{eqnarray*}
&&\varepsilon_k\cdot(B-\frac{2r+s}{2r+2s}) e_l\\\nonumber
&=&\sum_{1\leq i\leq r}\frac{r-i}{r}\mu^{i(l-k)}+\sum_{1\leq i\leq
r+s}\frac{i}{r+s}\mu^{i(l-k)}-\frac{2r+s}{2r+2s}\sum_{1\leq i\leq
r+s}\mu^{i(l-k)}.
\end{eqnarray*}
If $k=l$, $$\varepsilon_k\cdot(B-\frac{2r+s}{2r+2s}) e_l=\sum_{1\leq
i\leq r}\frac{r-i}{r}+\sum_{1\leq i\leq
r+s}\frac{i}{r+s}-\frac{2r+s}{2r+2s}\sum_{1\leq i\leq r+s}1=0.$$
Suppose $k\neq l$ and $r+s|(l-k)d$. Let $\xi=\mu^{l-k}$. Then
$\xi^d=1$ and $\xi\neq1$. For any $d|n$, we have $\sum_{1\leq i\leq
n}\xi^i=0$ and hence $\sum_{1\leq i\leq
n}i\xi^i=\frac{n}{d}\sum_{1\leq i\leq d}i\xi^i$. So we have
\begin{eqnarray*}
\varepsilon_k\cdot(B-\frac{2r+s}{2r+2s})
e_l&=&-\frac{1}{r}\sum_{1\leq i\leq
r}i\xi^i+\frac{1}{r+s}\sum_{1\leq i\leq
r+s}i\xi^i\\
&=&-\frac{1}{r}\frac{r}{d}\sum_{1\leq i\leq
d}i\xi^i+\frac{1}{r+s}\frac{r+s}{d}\sum_{1\leq i\leq d}i\xi^i=0.
\end{eqnarray*}
So $\varepsilon_k\cdot(B-\frac{2r+s}{2r+2s})v_0=0$ if $r+s|(k-j)d.$
Therefore $(B-\frac{2r+s}{2r+2s})v_0=(A_0-\alpha_0)v$ for some $v\in
V$. Then $v_0'=v_0,\ldots,v_{s-1}'=v_{s-1},v_s'=v_s-v$ satisfy the
system of equations (\ref{54}). Reversing the above argument, we get
the if part. So for the last assertion, it suffices to show that the
same assertion holds for the following system of equations
\begin{eqnarray}\label{eqn}
\sum_{0\leq i\leq k}(A_i-\alpha_i)v_{k-i}=0\hbox{ for any }0\leq
k\leq s.
\end{eqnarray}
Suppose the system of equations (\ref{eqn}) holds for some
$\alpha_0,\ldots,\alpha_s\in k$ and some $v_0,\ldots,v_s\in V$ with
$v_0\neq0$. There exists an integer $1\leq j\leq r+s$ such that
$\alpha_0=\mu^{rj}a_0^{\frac{r}{r+s}}$ and
$v_0=\sum_{r+s|(i-j)d}\gamma_ie_i$ for some $\gamma_i\in k$ with
$\gamma_j\neq0$. The system of equations (\ref{eqn}) is equivalent
to the following equation
$$\Big(\sum_{0\leq i\leq s}A_iZ'^i-\sum_{0\leq i\leq s}\alpha_iZ'^i\Big)\Big(\sum_{0\leq i\leq s}v_iZ'^i\Big)\equiv
0 \hbox { mod. } Z'^{s+1}.$$ There exist
$\rho_0=\mu^ja_0^{\frac{1}{r+s}},\rho_1,\ldots,\rho_s\in k$ such
that
$$\sum_{0\leq i\leq s}\alpha_iZ'^i\equiv\Big(\sum_{0\leq i\leq
s}\rho_iZ'^i\Big)^r\hbox { mod. } Z'^{s+1}.$$ Let
$$\Gamma=\left(\begin{array}{cccccc}
0&1&\\
\vdots&&\ddots\\
0&&&\ddots\\
a_sZ'^s&&&&\ddots\\
\vdots&&&&&1\\
 a_0&&&&&0
\end{array}\right)\hbox{ and }
\Gamma_0=\left(\begin{array}{cccc}
0&1&\\
\vdots&&\ddots\\
0&&&1\\
 a_0&&&0
\end{array}\right).$$
Then $\sum_{0\leq i\leq s}A_iZ'^i=\Gamma^r,\;A_0=\Gamma^r_0$ and
hence
$$\Big(\Gamma-\sum_{0\leq i\leq s}\rho_iZ'^i\Big)\Big(\sum_{0\leq k\leq r-1}\Big(\sum_{0\leq i\leq s}\rho_iZ'^i\Big)^k\Gamma^{r-1-k}\Big)
\Big(\sum_{0\leq i\leq s}v_iZ'^i\Big)\equiv0\hbox { mod. } Z'^{s+1}.
$$
Write
$$\Big(\sum_{0\leq k\leq r-1}\Big(\sum_{0\leq i\leq
s}\rho_iZ'^i\Big)^k\Gamma^{r-1-k}\Big) \Big(\sum_{0\leq i\leq
s}v_iZ'^i\Big)=\sum_{0\leq i}u_iZ'^i$$ for some $u_i\in V.$ Then
\begin{eqnarray*}
u_0&=&\sum_{0\leq k\leq
r-1}\rho_0^k\Gamma_0^{r-1-k}\sum_{r+s|(i-j)d}\gamma_ie_i\\
&=&\sum_{r+s|(i-j)d}\gamma_i\cdot\sum_{0\leq k\leq
r-1}\mu^{jk}a_0^{\frac{k}{r+s}}\mu^{i(r-1-k)}a_0^{\frac{r-1-k}{r+s}}e_i\\
&=&r\mu^{j(r-1)}\gamma_ja_0^{\frac{r-1}{r+s}}e_j\neq0.
\end{eqnarray*}
and
\begin{eqnarray}\label{M}
\Big(\Gamma-\sum_{0\leq i\leq s}\rho_iZ'^i\Big)\Big(\sum_{0\leq
i\leq s}u_iZ'^i\Big)\equiv 0\hbox{ mod. }Z'^{s+1}.
\end{eqnarray}
Since $\rho_0$ is a simple root of the characteristic polynomial of
$\Gamma_0$, by Hensel's lemma, $\rho_1,\ldots,\rho_s$ are uniquely
determined by $\rho_0$. So $\alpha_0,\ldots,\alpha_s$ are uniquely
determined by $\rho_0=\mu^ja_0^\frac{1}{r+s}$ ($1\leq j\leq r+s$).
This proves the last assertion.
\end{proof}

Now we are ready to prove Theorem 1$'$. By Remark \ref{r}, we assume
that $a(T)=\sum_{0\leq i\leq s}a_iT^i.$ Then the first equation of
(\ref{21}) means that $\frac{T}{Z'}$ is a root in $k[[Z']]$ of the
polynomial
$$\lambda^{r+s}-\sum_{0\leq i\leq s}a_iZ'^i\lambda^i\in k[[Z']][\lambda].$$ This polynomial is exactly
the characteristic polynomial of $\Gamma$. The characteristic
polynomial of $\Gamma_0$ is the polynomial $\lambda^{r+s}-a_0$ which
has no multiple roots, then by Hensel's lemma, $\Gamma$ has an
eigenvector $\sum_{i\geq 0}Z'^iv_i$ corresponding this eigenvalue
$\frac{T}{Z'}$ with $v_0\neq 0$. Since $\sum_{0\leq i\leq
s}Z'^iA_i=\Gamma^r$, we have
$$\Big(\sum_{0\leq i\leq s}Z'^iA_i\Big)\Big(\sum_{0\leq i} Z'^iv_i\Big)=\Big(\frac{T}{Z'}\Big)^r\Big(\sum_{0\leq i}Z'^iv_i\Big)
=\Big(\sum_{0\leq i}b_iZ'^i\Big)\Big(\sum_{0\leq i}Z'^iv_i\Big).$$
So
$$\sum_{0\leq i\leq k}(A_i-b_i)v_{k-i}=0 \hbox { for any } 0\leq k\leq s.$$
Recall that $\sum_{0\leq i\leq s}a_iT^{i-s}$ is assumed to be
irreducible. Then by Lemma \ref{l53} and \ref{l54}, we have
\begin{eqnarray*}
&&\sF^{(0,\infty)}\Big([r]_*[-r(a_0T^{-s}+a_1T^{1-s}+\ldots+a_s)]\Big)\\
&=&[r+s]_*[-(r+s)(b_0Z'^{-s}+b_1Z'^{1-s}+\ldots+b_s-\frac{2r+s}{2r+2s})]\\
&=&[r+s]_*[-(r+s)(b_0Z'^{-s}+b_1Z'^{1-s}+\ldots+b_s)+\frac{s}{2}].
\end{eqnarray*}

Suppose $r>s$. Given a formal Laurent series $\alpha$ in the
variable $\frac{1}{\sqrt[r]t}$ of order $-s$, consider the system of
equations (\ref{2}). We express $\frac{1}{\sqrt[r]t}$ as a formal
power series in $\sqrt[r-s]{t'}$ of order 1 using the first
equation, and then substitute this expression into the second
equation to get $\beta\in k((\sqrt[r-s]{t'}))$. Similar to equation
(\ref{*}), we have $\p_{t'}(\beta)=t.$ It follows that $\beta$ is a
formal Laurent series in $\sqrt[r-s]{t'}$ of order $-s$. Let
$Z=\frac{1}{\sqrt[r]t}$ and let $T'=\sqrt[r-s]{t'}$. Set
$$a(Z)=Z^st\p_t(\alpha) \hbox { and } b(T')=-T'^st'\p_{t'}\beta.$$ Then $a(Z)$ is
a formal power series in $Z$ of order 0 and $b(T')$ is a formal
power series in $T'$ of order 0. From the system of equations
(\ref{2}), we get
\begin{eqnarray}\label{22}\left\{\begin{array}{ll}
a(Z)=-(\frac{T'}{Z})^{r-s}\\
b(T')=-(\frac{T'}{Z})^r.
\end{array}\right.\end{eqnarray}
Similar to Theorem 1 and 1$'$, to prove Theorem 2, it suffices to
show the following theorem.
\bigskip

\noindent\textbf{Theorem 2$'$.} Suppose $r>s$. Given a formal power
series $a(Z)=\sum_{i\geq0}a_iZ^i$ with $a_i\in k$ and $a_0\neq0$,
suppose $b(T')=\sum_{i\geq0}b_iT'^i$ with $b_i\in k$ is a solution
of the system of equations (\ref{22}). We have
$b_s=\frac{r}{r-s}a_s$ and
\begin{eqnarray*}
&&\sF^{(\infty,0)}\Big([r]_*[-r(a_0Z^{-s}+a_1Z^{1-s}+\ldots+a_s)]\Big)\\
&=&[r-s]_*[-(r-s)(b_0T'^{-s}+b_1T'^{1-s}+\ldots+b_s)+\frac{s}{2}].
\end{eqnarray*}

\begin{proof}
The proof of $b_s=\frac{r}{r-s}a_s$ is similar to that of Theorem
$1'.$ Using the first equation of (\ref{22}), we can express $Z$ as
a formal power series in the variable $T'$ of order 1. We then
substitute this expression into the second equation to get $b(T')$
is a formal power series in $T'$ with nonzero constant term. That
is, $b_0\neq0$. Let $\zeta$ be an $r$-th root of $-1$ in $k$ and let
$Z=\zeta\cdot Z_1$. Let $[-]:k((z))\to k((z))$ be the automorphism
of $k$-algebra defined by $z\mapsto-z.$ From the system of equations
(\ref{22}), we get
\begin{eqnarray*}\left\{\begin{array}{ll}
\sum_{i\geq0}b_iT'^i=(\frac{T'}{Z_1})^r\\
\sum_{i\geq0}\zeta^{i-s}a_iZ_1^i=(\frac{T'}{Z_1})^{r-s}.
\end{array}\right.\end{eqnarray*}
Since $b_0\neq0$, by Theorem 1$'$, we have
\begin{eqnarray*}
&&\sF^{(0,\infty)}\Big([r-s]_*[-(r-s)(b_0T'^{-s}+b_1T'^{1-s}+\ldots+b_s)+\frac{s}{2}]\Big)\\
&=&[r]_*[-r(\zeta^{-s}a_0Z^{-s}+\zeta^{1-s}a_1Z^{1-s}+\ldots+a_s)+\frac{s}{2}+\frac{s}{2}]\\
&=&[-]^*[r]_*[-r(a_0Z^{-s}+a_1Z^{1-s}+\ldots+a_s)]\\
&=&\sF^{(0,\infty)}\Big(\sF^{(\infty,0)}\Big([r]_*[-r(a_0Z^{-s}+a_1Z^{1-s}+\ldots+a_s)]\Big)\Big).
\end{eqnarray*}
The theorem holds by \cite{BE}, Proposition 3.10.
\end{proof}

\section{Proof of Theorem 3}
Suppose $r<s$. Given a formal Laurent series $\alpha$ in the
variable $\frac{1}{\sqrt[r]t}$ of order $-s$, consider the system of
equations (\ref{2.}). We express $\frac{1}{\sqrt[r]t}$ as a formal
Laurent series in $\frac{1}{\sqrt[s-r]{t'}}$ of order 1 using the
first equation and then substitute this expression into the second
equation to get $\beta\in k((\frac{1}{\sqrt[s-r]{t'}}))$. Similar to
equation (\ref{*}), we have $\p_{t'}(\beta)=t.$ It follows that
$\beta$ is a formal Laurent series in $\frac{1}{\sqrt[s-r]{t'}}$ of
order $-s$. Let $Z=\frac{1}{\sqrt[r]t}$ and
$Z'=\frac{1}{\sqrt[s-r]{t'}}$. Set
$$a(Z)=Z^st\p_t(\alpha) \hbox{ and }
b(Z')=Z'^st'\p_{t'}(\beta).$$ Then $a(Z)$ is a formal power series
in $Z$ of order 0 and $b(Z')$ is a formal power series in $Z'$ of
order 0. From the system of equations (\ref{2.}), we get
\begin{eqnarray}\label{23}\left\{\begin{array}{ll}
a(Z)=-(\frac{Z}{Z'})^{s-r}\\
b(Z')=(\frac{Z'}{Z})^r.
\end{array}\right.\end{eqnarray}
Similar to Theorem 1 and 1$'$, to prove Theorem 3, it suffices to
show the following theorem.

\noindent\textbf{Theorem 3$'$.} Suppose $s>r$. Given a formal power
series $a(Z)=\sum_{i\geq0}a_iZ^i$ with $a_i\in k$ and $a_0\neq 0$,
solve the system of equations (\ref{23}) to get
$b(Z')=\sum_{i\geq0}b_iZ'^i$ for some $b_i\in k$. Then
$b_s=\frac{r}{s-r}a_s$ and
\begin{eqnarray*}
&&\sF^{(\infty,\infty)}\Big([r]_*[-r(a_0Z^{-s}+a_1Z^{1-s}+\ldots+a_s)]\Big)\\
&=&[s-r]_*[-(s-r)(b_0Z'^{-s}+b_1Z'^{1-s}+\ldots+b_s)+\frac{s}{2}].
\end{eqnarray*}

\begin{lem}\label{l58}
Set $h=a_0Z^{-s}+a_1Z^{1-s}+\ldots+a_s$. We can reduce Theorem 3$'$
to the case where $s\geq2r$ and where $h$ is irreducible with
respect to the Galois extension $k((Z))/k((z))$.
\end{lem}

\begin{proof}
The proof of $b_s=\frac{r}{s-r}a_s$ is similar to that of Theorem
1$'$ and the proof of the last assertion is similar to that of Lemma
\ref{l52}. If $s<2r$, then $s>2(s-r)$. Let $\zeta$ be an $r$-th root
of $-1$ in $k$ and let $Z=\zeta\cdot Z_1$. From the system of
equations (\ref{23}), we get
\begin{eqnarray*}\left\{\begin{array}{ll}
\sum_{i\geq0}b_iZ'^i=-(\frac{Z'}{Z_1})^r\\
\sum_{i\geq0}\zeta^{i-s}a_iZ_1^i=(\frac{Z_1}{Z'})^{s-r}.
\end{array}\right.\end{eqnarray*}
We prove $b_0\neq0$ similarly as in Theorem 2$'$. Applying this
theorem to $[s-r]_*[-(s-r)(b_0Z'^{-s}+\ldots+b_s)+\frac{s}{2}]$, we
have
\begin{eqnarray*}
&&\sF^{(\infty,\infty)}\Big([s-r]_*[-(s-r)(b_0Z'^{-s}+b_1Z'^{1-s}+\ldots+b_s)+\frac{s}{2}]\Big)\\
&=&[r]_*[-r(\zeta^{-s}a_0Z^{-s}+\zeta^{1-s}a_1Z^{1-s}+\ldots+a_s)+\frac{s}{2}+\frac{s}{2}]\\
&=&[-]^*[r]_*[-r(a_0Z^{-s}+a_1Z^{1-s}+\ldots+a_s)]\\
&=&\sF^{(\infty,\infty)}\Big(\sF^{(\infty,\infty)}\big([r]_*[-rh]\big)\Big).
\end{eqnarray*}
The lemma holds by \cite{BE}, Proposition 3.12 (iv).
\end{proof}
From now on, we assume $h$ is irreducible.

\medskip

Let's describe the formal connection
$\sF^{(\infty,\infty)}\big([r]_*[-rh]\big)$ on $k((z')).$

The formal connection $[-rh]$ on $k((Z))$ consist of a one
dimensional $k((Z))$-vector space with a basis $e'$ and a $k$-linear
map $Z\p_Z:k((Z))e'\to k((Z))e'$ satisfying
$$Z\p_Z(ge')=(Z\p_Z(g)-rhg)e'$$
for any $g\in k((Z)).$ Since the formal connection $[-rh]$ on
$k((Z))$ has slope $s$, we get $k[[Z]] e',\;Z^{-s}k[[Z]] e'$ is a
good lattices pair for it. Identify $[r]_*[-rh]$ with $k((Z)) e'$ as
$k((z))$-vector spaces. So the connection $[r]_*[-rh]$ on $k((z))$
has pure slope $\frac{s}{r}$ and $k[[Z]] e',\;Z^{-s}k[[Z]] e'$ is a
good lattices pair for this connection. The action of the
differential operator $z\p_z$ on $k((Z))e'$ is given by
$$z\p_z(ge')=(z\p_z(g)-hg)e'$$
for any $g\in k((Z)).$ Then for any $i\in\Z$, we have
$$z^2\p_z(Z^{-(r+i)}e')=-\frac{r+i}{r}Z^{-i}e'-(a_0Z^{-(i+s)}e'+\ldots+a_sZ^{-i}e').$$
By \cite{BE}, Proposition 3.12 (ii), the map
$$\iota:k((Z))e'\to\sF^{(\infty,\infty)}\Big([r]_*[-rh]\Big)$$ is an
isomorphism of $k$-vector spaces. As in \cite{BE}, Proposition 3.14,
$(\iota Z^{-1}e',\ldots,\iota Z^{-(s-r)}e')$ is a basis of
$\sF^{(\infty,\infty)}\Big([r]_*[-rh]\Big)$ over $k((Z'))$. By the
relation $\iota\circ z^2\p_z=\frac{1}{z'}\circ\iota$ and
$-\iota\circ\frac{1}{z}=z'^2\p_{z'}\circ\iota$ in \cite{BE},
Proposition 3.12 (iii), we have
\begin{eqnarray*}
z'^2\p_{z'}(\iota Z^{-(i+s-r)}e')&=&-\iota
Z^{-(i+s)}e'\\&=&\frac{a_s}{a_0}\iota
Z^{-i}e'+\ldots+\frac{a_1}{a_0}\iota Z^{-(i+s-1)}e'
+\frac{1}{a_0z'}\iota Z^{-(r+i)}e'+\frac{r+i}{ra_0}\iota Z^{-i}e'.
\end{eqnarray*}
Let
$$A=\left(\begin{array}{ccccccc}
0&&&&&&-\frac{a_s}{a_0}\\
1&&&&&&\vdots\\
&\ddots&&&&&-\frac{a_{s-r+1}}{a_0}\\
&&\ddots&&&&-\frac{a_{s-r}}{a_0}-\frac{1}{a_0z'}\\
&&&\ddots&&&-\frac{a_{s-r-1}}{a_0}\\
&&&&\ddots&&\vdots\\
&&&&&1&-\frac{a_1}{a_0}
\end{array}\right).$$
For any $i\in\Z$, let $B_i$ be the $s\times s$-matrix whose entries
are all zero except the $(1, s)$-th entry which is valued by
$-\frac{r+i}{ra_0}$. We have
$$(\iota Z^{-(i+1)}e',\ldots,\iota Z^{-(i+s)}e')=(\iota Z^{-i}e',\ldots,\iota Z^{-(i+s-1)}e')(A+B_i).$$
So
\begin{eqnarray*}z'^2\p_{z'}(\iota Z^{-1}e',\ldots,\iota
Z^{-s}e')&=&-(\iota Z^{-(r+1)}e',\ldots,\iota
Z^{-(r+s)}e')\\&=&-(\iota Z^{-1}e',\ldots,\iota
Z^{-s}e')\prod_{1\leq i\leq r}(A+B_i).
\end{eqnarray*}
 Consider the connection
$$[s-r]^*\Big(\sF^{(\infty,\infty)}\Big([r]_*[-rh]\Big)\Big)=k((Z'))\otimes_{k((z'))}\sF^{(\infty,\infty)}\Big([r]_*[-rh]\Big).$$
Set $\wedge={\rm diag}\{Z',\ldots,Z'^s\}$ and
$\varepsilon'=(Z'\otimes\iota Z^{-1}e', \ldots,Z'^s\otimes\iota
Z^{-s}e')$. We have
\begin{eqnarray*}
Z'\p_{Z'}(\varepsilon')&=&\varepsilon'\cdot\Big({\rm
diag}\{1,\ldots,s\}-\frac{s-r}{z'}\wedge^{-1}\Big(\prod_{1\leq
i\leq r}(A+B_i)\Big)\wedge\Big)\\
&=&\varepsilon'\cdot\Big({\rm
diag}\{1,\ldots,s\}-\frac{s-r}{Z'^s}\prod_{1\leq i\leq
r}\Big(Z'\wedge^{-1}(A+B_i)\wedge\Big)\Big).
\end{eqnarray*}
We have
$$Z'\wedge^{-1}A\wedge=\left(\begin{array}{ccccccc}
0&&&&&&-\frac{a_s}{a_0}Z'^s\\
1&&&&&&\vdots\\
&\ddots&&&&&-\frac{a_{s-r+1}}{a_0}Z'^{s-r+1}\\
&&\ddots&&&&-\frac{a_{s-r}}{a_0}Z'^{s-r}-\frac{1}{a_0}\\
&&&\ddots&&&-\frac{a_{s-r-1}}{a_0}Z'^{s-r-1}\\
&&&&\ddots&&\vdots\\
&&&&&1&-\frac{a_1}{a_0}Z'
\end{array}\right)$$
and $Z'\wedge^{-1}B_i\wedge$ is the $s\times s$-matrix whose entries
are all zero except the $(1,s)$-th entry which is valued by
$-\frac{r+i}{ra_0}Z'^s$. So we can write
$${\rm diag}\{1,2,\ldots,s\}-\frac{s-r}{Z'^s}\prod_{1\leq i\leq
r}\Big(Z'\wedge^{-1}(A+B_i)\wedge\Big)=-(s-r)\sum_{i\geq0}Z'^{i-s}C_i$$
and $$\big(Z'\wedge^{-1}A\wedge\big)^r=\sum_{i\geq0}Z'^iC'_i$$ for
some matrices $C_i$ and $C'_i$ with entries in $k$. Then $C_i=C'_i$
for all $0\leq i\leq s-1$ and $$C'_s-C_s={\rm
diag}\{\frac{1}{s-r},\ldots,\frac{s}{s-r}\}-P$$ where $P$ is the
$s\times s$-matrix whose entries are all zero except the
$(i,i+s-r)$-th entry which is valued by $-\frac{r+i}{ra_0}$ $(1\leq
i\leq r).$ Let $W$ be the $k$-vector space of column vectors in $k$
of length $s.$ We have

\begin{lem}\label{l5.9}
Suppose $s\geq 2r$ and $h$ is irreducible with respect to the Galois
extension $k((Z))/k((z))$. Given $\alpha_0,\ldots,\alpha_s\in k$
with $\alpha_0\neq0$, the following three conditions are equivalent:

$(1)$
$\sF^{(\infty,\infty)}\big([r]_*[-rh]\big)=[s-r]_*[-(s-r)\sum_{0\leq
i\leq s}\alpha_iZ'^{i-s}].$

$(2)$ $[-(s-r)\sum_{0\leq i\leq s}\alpha_iZ'^{i-s}]$ is a
subconnection of $[s-r]^*\sF^{(\infty,\infty)}\big([r]_*[-rh]\big).$

$(3)$ There exist $N\in\Z$ and $w_0,\ldots,w_s\in W$ such that
$w_0\neq0$ and
\begin{eqnarray}\left\{\begin{array}{ll}
\sum_{0\leq i\leq k}(C_i-\alpha_i)w_{k-i}=0\;(0\leq k\leq s-1),\\
\sum_{0\leq i\leq
s-1}(C_i-\alpha_i)w_{s-i}+(C_s-\alpha_s-\frac{N}{s-r})w_0=0.
\end{array}\right.\end{eqnarray}
\end{lem}

\begin{proof}
Set $U=W\otimes_kk((Z'))$ and $\sW=W\otimes_kk[[Z']]$. Let
$u=(u_1,\ldots,u_s)$ be the canonical basis of $W.$ There exists a
unique connection $(U,Z'\p_{Z'})$ such that the action of
$Z'\p_{Z'}$ on elements of $W$ can be written as
$$Z'\p_{Z'}(w)=-(s-r)\sum_{i\geq 0}Z'^{i-s}C_i(w).$$
The map of $k((Z'))$-vector spaces
$$U\to[s-r]^*\Big(\sF^{(\infty,\infty)}\Big([r]_*[-rh]\Big)\Big)$$ which maps each $u_i$
to $Z'^i\otimes\iota Z^{-i}e'$ is a surjective morphism of
connections. We have $Z'^{s+1}\p_{Z'}(\sW)\subset\sW$. Let
$\psi:\sW\to\sW/Z'\sW\cong W$ be the canonical map. The $k$-linear
action on $W\cong\sW/Z'\sW$ induced by $Z'^{s+1}\p_{Z'}$ is
$-(s-r)C_0.$ Write
$$Z'\wedge^{-1}A\wedge=\sum_{i\geq0}Z'^iD_i$$ for some
matrices $D_i$ with entries in $k$. The characteristic polynomial of
$D_0$ is $\lambda^s+\frac{1}{a_0}\lambda^r$. So $W$ is the direct
sum of two subspaces $W_0$ and $W_1$, invariant under $D_0,$ and
such that $D_0|_{W_0}$ is nilpotent, $D_0|_{W_1}$ is invertible.
Then dim$W_0=r$ and dim$W_1=s-r$. Since $C_0=D_0^r$, we have $W_0$
and $W_1$ are $C_0$-invariant, and then $C_0|_{W_0}=0$, $C_0|_{W_1}$
is invertible. By the splitting lemma in \cite{Le}, 2, $\sW$ is the
direct sum of two free submodules $\sW_0$ and $\sW_1$, invariant
under $Z'^{s+1}\p_{Z'}$, and such that
$W_0=\psi(\sW_0),\;W_1=\psi(\sW_1).$ Let $U_0$, $U_1$ be the
subconnections of $U$ generated by $\sW_0$, $\sW_1$, respectively.
Then $U=U_0\oplus U_1.$ The induced action of $Z'^{s+1}\p_{Z'}$ on
$W_0$ is 0, so the slopes of the connection $U_0$ are all $<s$. But
$[s-r]^*\Big(\sF^{(\infty,\infty)}\Big([r]_*[-rh]\Big)\Big)$ is an
$s-r$ dimensional connection on $k((Z'))$ with pure slope $s$, we
have
$${\rm Hom}_{{\rm
conn.}}\Big(U_0,[s-r]^*\Big(\sF^{(\infty,\infty)}\Big([r]_*[-rh]\Big)\Big)\Big)=(0)$$
and then
$$U_1\cong[s-r]^*\Big(\sF^{(\infty,\infty)}\Big([r]_*[-rh]\Big)\Big).$$ For any
one dimensional formal connection $L$ on $k((Z'))$ with slope $s$,
we have $${\rm Hom}_{\rm conn.}(L,U_0)=(0)$$ and then
\begin{eqnarray*}
{\rm Hom}_{\rm conn.}(L,U)&=&{\rm Hom}_{\rm
conn.}(L,U_1)\bigoplus{\rm Hom}_{\rm conn.}(L,U_0)\\&=&{\rm
Hom}_{\rm
conn.}\Big(L,[s-r]^*\Big(\sF^{(\infty,\infty)}\Big([r]_*[-rh]\Big)\Big)\Big).
\end{eqnarray*}
 So to find a one dimensional subconnection in
$[s-r]^*\Big(\sF^{(\infty,\infty)}\Big([r]_*[-rh]\Big)\Big)$ is
equivalent to finding a one dimensional  subconnection in $U$ of
slope $s$. By Lemma \ref{l5.10}, the remainder proof is similar to
that of Lemma \ref{l53}.
\end{proof}

\begin{lem}\label{l5.10}
Suppose $s\geq2r$. Given $\alpha_0,\ldots,\alpha_s\in k$ with
$\alpha_0\neq0$, there exist $w_0,\ldots,w_s\in W$ such that
$w_0\neq0$ and
\begin{eqnarray}\label{5.10}
\sum_{0\leq i\leq k}(C_i-\alpha_i)w_{k-i}=0\;(0\leq k\leq s)
\end{eqnarray}
if and only if there exist $w_0',\ldots,w_s'\in W$ such that
$w_0'\neq0$ and
\begin{eqnarray}\label{5.11}\left\{\begin{array}{ll}
\sum_{0\leq i\leq k}(C'_i-\alpha_i)w'_{k-i}=0\;(0\leq k\leq s-1),\\
\sum_{0\leq i\leq
s-1}(C'_i-\alpha_i)w'_{s-i}+(C'_s-\alpha_s-\frac{s-2r}{2s-2r})w_0'=0.
\end{array}\right.\end{eqnarray}
Moreover, there are at most $s-r$ $(s+1)$-tuples
$(\alpha_0,\ldots,\alpha_s)$ in $k$ such that $\alpha_0\neq0$ and
the system of equations (\ref{5.10}) (resp. (\ref{5.11})) holds for
some $w_0,\ldots,w_s\in W$ with $w_0\neq 0$ (resp.
$w_0',\ldots,w'_s\in W$ with $w'_0\neq0$).
\end{lem}

\begin{proof}
Let $\eta$ be a primitive $(s-r)$-th root of unity in $k$. We fix an
$(s-r)$-th root $(-a_0)^{\frac{1}{s-r}}$ of $-a_0$. For any $1\leq
j\leq s-r$, set $e'_j$ to be the column vector
$(0,\ldots,0,\eta^{(r+1)j}(-a_0)^{\frac{r+1}{s-r}},\ldots,\eta^{sj}(-a_0)^{\frac{s}{s-r}})$
and $\varepsilon_j'$ the row vector
$(\eta^{-j}(-a_0)^{-\frac{1}{s-r}},\ldots,\eta^{-sj}(-a_0)^{-\frac{s}{s-r}})$.
We have
$$C_0\cdot e'_j=\eta^{-rj}(-a_0)^{-\frac{r}{s-r}}\cdot e'_j,\; \varepsilon'_j\cdot C_0=
\eta^{-rj}(-a_0)^{-\frac{r}{s-r}}\cdot\varepsilon'_j,\;
\varepsilon'_i\cdot e'_j=(s-r)\delta_{ij}.$$ Set $d=(r,s)$. Let
$W_0$ be as in Lemma \ref{l5.9}. We have $C_0|_{W_0}=0.$ Then ${\rm
ker}(C_0-\eta^{-rj}(-a_0)^{-\frac{r}{s-r}})$ is generated by those
$e'_k$ with $s-r|(k-j)d$, and ${\rm
im}(C_0-\eta^{-rj}(-a_0)^{-\frac{r}{s-r}})$ is generated by $W_0$
and the other $e'_k$'s.  So
$${\rm im}(C_0-\eta^{-rj}(-a_0)^{-\frac{r}{s-r}})=\{w\in W|\varepsilon'_k\cdot w=0\hbox{ for all } k \hbox{
satisfying }s-r|(k-j)d\}.$$ For the only if part, suppose the system
of equations (\ref{5.10}) holds for some $w_0,\ldots,w_s\in W$ with
$w_0\neq0$. So $\alpha_0=\eta^{-rj}(-a_0)^{-\frac{r}{s-r}}$ for some
integer $j$ and then $w_0=\sum_{s-r|(i-j)d}\sigma_ie'_i$ for some
$\sigma_i\in k$. Since $s\geq2r$, for any $1\leq k,\;l\leq s-r$, we
have
\begin{eqnarray*}
&&\varepsilon'_k\cdot({\rm
diag}\{\frac{1}{s-r},\ldots,\frac{s}{s-r}\}-P-\frac{s-2r}{2s-2r})e'_l\\
&=&\sum_{r+1\leq i\leq s}\frac{i}{s-r}\eta^{(l-k)i}-\sum_{1\leq
i\leq r}\frac{r+i}{r}\eta^{(l-k)i}-\frac{s-2r}{2}\delta_{kl}.
\end{eqnarray*}
If $k=l$, then
\begin{eqnarray*}
&&\varepsilon'_k\cdot({\rm
diag}\{\frac{1}{s-r},\ldots,\frac{s}{s-r}\}-P-\frac{s-2r}{2s-2r})e'_l\\
&=&\sum_{r+1\leq i\leq s}\frac{i}{s-r}-\sum_{1\leq i\leq
r}\frac{r+i}{r}-\frac{s-2r}{2}=0.
\end{eqnarray*}
If $k\neq l$ and $s-r|(l-k)d$, then $(\eta^{l-k})^d=1$ and
$\eta^{l-k}\neq1$. We have
\begin{eqnarray*}
&&\varepsilon'_k\cdot({\rm
diag}\{\frac{1}{s-r},\ldots,\frac{s}{s-r}\}-P-\frac{s-2r}{2s-2r})e'_l\\
&=&\frac{s-r}{d}\sum_{1\leq i\leq
d}\frac{i}{s-r}\eta^{(l-k)i}-\frac{r}{d}\sum_{1\leq i\leq
d}\frac{i}{r}\eta^{(l-k)i}=0.
\end{eqnarray*}
 So $\varepsilon'_k\cdot({\rm
diag}\{\frac{1}{s-r},\ldots,\frac{s}{s-r}\}-P-\frac{s-2r}{2s-2r})w_0=0$
if $s-r|(k-j)d$. Therefore $$({\rm
diag}\{\frac{1}{s-r},\ldots,\frac{s}{s-r}\}-P-\frac{s-2r}{2s-2r})w_0=(C_0-\alpha_0)w$$
for some $w\in W$. Then
$w_0'=w_0,\ldots,w'_{s-1}=w_{s-1},w_s'=w_s-w$ satisfy the system of
equations (\ref{5.11}). Reversing the above argument, we get the if
part. The characteristic polynomial of $D_0$ is
$\lambda^s+\frac{1}{a_0}\lambda^r$. Each nonzero root of this
polynomial is simple. Since
$\sum_{i\geq0}Z'^iC'_i=(Z'\wedge^{-1}A\wedge)^r$ and $C_0=D_0^r$,
the proof of the last assertion is similar to that of Lemma
\ref{l54}.
\end{proof}

Now we are ready to prove Theorem 3$'$. Similar to Remark \ref{r},
we assume $a(Z)=\sum_{0\leq i\leq s}a_iZ^i$. Then the first equation
of (\ref{23}) means that $\frac{Z'}{Z}$ is a root in $k[[Z']]$ with
nonzero constant term of the polynomial
$$\lambda^s+\frac{a_1}{a_0}Z'\lambda^{s-1}+\ldots+\frac{a_s}{a_0}Z'^s+\frac{1}{a_0}\lambda^r\in k[[Z']][\lambda].$$
This polynomial is exactly the characteristic polynomial of
$Z'\wedge^{-1}A\wedge.$ The characteristic polynomial of $D_0$ is
$\lambda^s+\frac{1}{a_0}\lambda^r$ which has no nonzero multiple
roots, by Hensel's lemma, $Z'\wedge^{-1}A\wedge$ has an eigenvector
$\sum_{i\geq0}Z'^iw_i$ corresponding this eigenvalue $\frac{Z'}{Z}$
with $w_0\neq0$. Since
$\sum_{i\geq0}Z'^iC'_i=(Z'\wedge^{-1}A\wedge)^r,$ we have
\begin{eqnarray*}
\Big(\sum_{i\geq0}Z'^iC'_i\Big)\Big(\sum_{i\geq0}Z'^iw_i\Big)
=\Big(\frac{Z'}{Z}\Big)^r\Big(\sum_{i\geq0}Z'^iw_i\Big)=\Big(\sum_{i\geq0}b_iZ'^i\Big)\Big(\sum_{i\geq0}Z'^iw_i\Big).
\end{eqnarray*}
That is,
$$\sum_{0\leq i\leq k}(C'_i-b_i)w_{k-i}=0\hbox{ for any
}k\geq0.$$ Recall that $s\geq2r$ and $\sum_{0\leq i\leq
s}a_iZ^{i-s}$ is assumed to be irreducible. By Lemma \ref{l5.9} and
\ref{l5.10}, we have
\begin{eqnarray*}
&&\sF^{(\infty,\infty)}\Big([r]_*[-rh]\Big)\\
&=&[s-r]_*[-(s-r)(b_0Z'^{-s}+b_1Z'^{1-s}+\ldots+b_s-\frac{s-2r}{2s-2r})]\\
&=&[s-r]_*[-(s-r)(b_0Z'^{-s}+b_1Z'^{1-s}+\ldots+b_s)+\frac{s}{2}].
\end{eqnarray*}
\newpage

\bibliographystyle{plain}


\begin{thebibliography}{99}
\bibitem{BBE} Beilinson, A.; Bloch, S.; Esnault, H.: $\epsilon$-factors
for Gau{\ss}-Manin determinants, Moscow Mathematical Journal, vol.
2, {\bf 3} (2002), 477-532.
\bibitem{BE} Bloch, S. and Esnault, H.: Local Fourier Transforms and
rigidity for $\mathcal D$-Modules, Asian J. Math. 8 (2004), no. 4,
587--605.
\bibitem{De} Deligne, P.: \'Equations Diff\'erentielles \`a Points Singuliers
R\'eguliers, Lectures Notes {\bf 163}, Springer Verlag.
\bibitem{F} Fu, Lei.: Calculation of $\ell$-adic local Fourier
transformations, arXiv:math/0702436.
\bibitem{Ka-1} Katz, N.: On the calculation of some differential Galois
groups, Inv. Math. {\bf 87} (1986), 13-61.
\bibitem{La} Laumon, G.: Transformation de Fourier, constantes d'\'equations
fonctionnelles et conjecture de Weil, Publ. Math. IHES {\bf 65}
(1987), 131-210.
\bibitem{Le} Levelt, A.H.M.: Jordan decomposition for a class of singualr differential operators,
Ark. Math. 13 (1975), 1-27.
\bibitem{Sa} Sabbah, Claude.:  An explicit stationary phase
formula for the local formal Fourier-Laplace transform,
arXiv:0706.3570.
\bibitem{S} Serre, J-P.: Local Fields, Graduate Texts in
Mathematics, Springer Verlag.


\end{thebibliography}
\end{document}